\documentclass[12pt]{article}
\usepackage{amsmath,amssymb}
\newtheorem{theorem}{Theorem}
\begin{document}
\begin{center}
\large
\textbf{Ergodic Density Estimates for some diffusion processes}
\end{center}
\begin{center}
\textbf{Bert Koehler and Volker Krafft}
\end{center}
\begin{abstract}
\noindent
For n-dimensional ergodic diffusion processes with values in $G=\mathbb{R}_{+}^n$ we prove time-independent upper bounds
for the transitional density and so also for the unique ergodic density. We do not require geodesic completeness
of the elliptic symbol towards the boundary of $G$. 
\end{abstract}
Let $W_{1,t},...,W_{n,t}$ be independent Brownian motions, let $G=\mathbb{R}_{+}^n$
and let $\sigma_{ij},\mu_i:G\longrightarrow\mathbb{R}$ be smooth functions. Assume there is a strong
solution process $X_t:\Omega\longrightarrow G$ of the equation
\begin{eqnarray*}
X_{i,t}=X_{i,0}+\int_0^t \mu_i(X_s)ds+\sum_{j=1}^n\int_0^t\sigma_{ij}(X_s)dW_{j,s}
\end{eqnarray*}
The image domain of nonnegative values appears naturally in interest rate or credit intensity modelling.
The generator of the process $X_t$ is given by
\begin{eqnarray*}
L(f)(x)=\sum_{i,j=1}^n a_{ij}(x_1,...,x_n)\frac{\partial^2 f}{\partial x_i\partial x_j}+
\sum_{i=1}^n\mu_i(x_1,...,x_n)\frac{\partial f}{\partial x_i}
\end{eqnarray*} 
where
\begin{eqnarray*}
a_{ij}(x)=\frac{1}{2}\sum_{k=1}^n\sigma_{ik}(x)\sigma_{jk}(x)=a_{ji}(x)
\end{eqnarray*}
We require
\begin{eqnarray*}
(a_{ij}(x))>0\mbox{ positive definite for all }x\in G
\end{eqnarray*}
Examples include so called affine processes where $\mu(x),a_{ij}(x)$ are affin linear functions such as
\begin{eqnarray*}
X_{i,t}=X_{i,0}+\int_0^t(\mu_{i0}+\sum_{j\neq i}\mu_{ij}X_{j,s}-\mu_{ii}X_{i,s})ds+\sigma_i\int_0^t\sqrt{X_{i,s}}dW_{i,s}
\end{eqnarray*}
but also more complicated processes with stochastic volatility such as
\begin{eqnarray*}
X_{1,t}=X_{1,0}+\int_0^t(\mu_{1,0}-\mu_{1,1}X_{1,s})ds+\sigma_1\int_0^t\sqrt{X_{1,s}}dW_{1,s}\\
X_{2,t}=X_{2,0}+\int_0^t(\mu_{2,0}+\mu_{2,1}X_{1,s}-\mu_{2,2}X_{2,s})ds+\sigma_2\int_0^t\sqrt{X_{1,s}X_{2,s}}dW_{2,s}\\
...\\
X_{n,t}=X_{n,0}+\int_0^t(\mu_{n,0}+\sum_1^{n-1}\mu_{n,j}X_{j,s}-\mu_{n,n}X_{n,s})ds+\sigma_n\int_0^t\sqrt{X_{n-1,s}X_{n,s}}dW_{n,s}
\end{eqnarray*}
Both classes will meet the conditions we will later impose on $\mu_i(x),a_{ij}(x)$ and in both cases the metric which corresponds to the
elliptic symbol $g_{ij}(x)=(a_{ij}(x))^{-1}$ is not geodesically complete towards the 0-boundary of $G$, so standard heat kernel estimates
do not apply directly.\\
\\
From pointwise positivity of the elliptic symbol $(a_{ij}(x))>0$ H\"ormanders criterion implies the existence of a smooth transitional density
$\rho(t,X_0,y)$ with respect to Lebesgue measure. The adjoint operator of $L$ with respect to euclidean metric is given by
\begin{eqnarray*}
L^{*}=\sum_{i,j=1}^n a_{ij}(y)\partial_i\partial_j-\sum_{j=1}^n b_j(y)\partial_j+c(y)
\end{eqnarray*}
where
\begin{eqnarray*}
b_j(y)=\mu_j(y)-2\sum_{i=1}^n\partial_i a_{ij}(y)\\
c(y)=\sum_{i,j=1}^n \partial_i\partial_j a_{ij}(y)-\sum_{j=1}^n\partial_j \mu_j(y)
\end{eqnarray*}
In the following we denote by $\partial_j=\frac{\partial}{\partial y_j}$ a partial derivative in the
state variables and $\partial_t=\frac{\partial}{\partial t}$.
The transition density fulfills Kolmogorovs equation
\begin{eqnarray*}
\partial_t\rho=L_y^{*}(\rho)
\end{eqnarray*}
For proving estimates of the transition density we impose the following conditions:\\
\\
(1) Ergodicity of the process $X_t$ is usually ensured by means of a barrier function. So let
$\psi:G\longrightarrow\mathbb{R}^{>0}$ be a smooth function with
\begin{eqnarray*}
\lim_{x_j\rightarrow 0}\psi(x_1,...,x_n)=0=\lim_{x_j\rightarrow\infty}\psi(x_1,...,x_n)\mbox{ for all }j=1,...,n
\end{eqnarray*}
and
\begin{eqnarray*}
L^{*}(\psi)(x)<-\psi(x)
\end{eqnarray*}
outside a compact set $K\subset\subset G$.
This is usually required for a barrier function to reflect a mean reversion property of the generator.\\
\\
(2) We need a further technical condition: There is a compact cube $K\subset\subset G$ such that for all
$y\in G\backslash K$
\begin{eqnarray*}
-\sum_{i,j=1}^n \partial_i\partial_j a_{ij}+\sum_{i=1}^n\partial_i\mu_i-
2\sum_{i,j=1}^n \partial_i(a_{ij}\partial_j\log\psi)\leq 0
\end{eqnarray*}
\\
(3) Condition (1) reads for $y\in G\backslash K$
\begin{eqnarray*}
\sum_{i,j}^n a_{ij}\partial_i\partial_j\psi-\sum_{j=1}^n\Big(\mu_j-2\sum_{i=1}^n\partial_i a_{ij}\Big)\cdot\partial_j\psi+
\Big(\sum_{i,j}^n\partial_i\partial_j a_{ij}-\sum_{j=1}^n\partial_j\mu_j\Big)\cdot\psi\leq-\psi
\end{eqnarray*}
Condition (2) yields for $y\in G\backslash K$
\begin{eqnarray*}
-2\sum_{i,j}^n a_{ij}\partial_i\partial_j\psi+\frac{2}{\psi}\cdot\sum_{i,j}^n a_{ij}\partial_i\psi\partial_j\psi-
2\sum_{i,j}^n\partial_i a_{ij}\partial_j\psi+
\Big(\sum_{j=1}^n\partial_j\mu_j-\sum_{i,j}^n\partial_i\partial_j a_{ij}\Big)\cdot\psi\leq 0
\end{eqnarray*}
Adding both inequalities we get outside a compact set $K\subset\subset G$
\begin{eqnarray*}
-\sum_{i,j}^n a_{ij}\partial_i\partial_j\psi+\frac{2}{\psi}\cdot\sum_{i,j}^n a_{ij}\partial_i\psi\partial_j\psi-
\sum_{j=1}^n\mu_j\partial_j\psi\leq-\psi
\end{eqnarray*}
But this is the same as
\begin{eqnarray*}
L\Big(\frac{1}{\psi}\Big)\leq-\frac{1}{\psi}+C\mbox{ on all of }G
\end{eqnarray*}
So if we define a function
\begin{eqnarray*}
F(t,X_0)=E\left(\frac{1}{\psi(X_t)}|X_0\right)
\end{eqnarray*}
(assuming for a moment that it exists) we get from Ito lemma
\begin{eqnarray*}
F(t,X_0)=\frac{1}{\psi(X_0)}+\int_0^t E\Big(L\Big(\frac{1}{\psi(X_s)}\Big)|X_0\Big)ds
\end{eqnarray*} 
and so
\begin{eqnarray*}
\frac{\partial F}{\partial t}=E\Big(L\Big(\frac{1}{\psi(X_t)}\Big)|X_0\Big)\leq-F(t,X_0)+C
\end{eqnarray*}
Now apply Gronwall inequality to get
\begin{eqnarray*}
F(t,X_0)\leq e^{-t}\cdot\left(\frac{1}{\psi(X_0)}+\int_0^t e^{s}C_R ds\right)<e^{-t}\cdot\frac{1}{\psi(X_0)}+C
\end{eqnarray*}
which shows that $F(t,X_0)$ is uniformly bounded for all times $t\geq 0$ (dependent on $X_0$).
Using the stopped version of Ito (called Dynkin lemma) one can now show that $F(t,X_0)$ exists and
is (depending on $X_0$) uniformly bounded in $t$.\\
\\
In fact we will also need that expressions like
\begin{eqnarray*}
E\left(\frac{1}{\psi(X_t)}\cdot\left(1+\sum_{i,j}^n a_{ij}(X_t)(1+\frac{1}{X_{it}})(1+\frac{1}{X_{jt}})\right)|X_0\right)<C_1
\end{eqnarray*}
are $t$-uniformly bounded which means that conditions (1),(2) should additionally be fulfilled by 
$\widetilde{\psi}(x)=x_i x_j\psi(x)$, so the asymptotics of $\psi(x)$ towards the 0-boundaries of $G$ has to be chosen
some levels below optimal parameters. In the following we will refer to this as assumption (3).\\
\\
For the class of affine processes one can choose
\begin{eqnarray*}
\psi(x)=x_1^{\beta_1}x_2^{\beta_2}...x_n^{\beta_n}\cdot e^{-\gamma_1 x_1-...-\gamma_n x_n}
\end{eqnarray*}
with suitable positive parameters $\beta_j,\gamma_j>0$ if $2\mu_{j,0}>\sigma_j^2$ as a barrier function.\\
For the second example a barrier function is given by
\begin{eqnarray*}
\psi(x)=x_1^{\beta_1}x_2^{\beta_2}...x_n^{\beta_n}\cdot e^{-\phi(x)}
\end{eqnarray*}
with
\begin{eqnarray*}
\phi(x)=\left(1+\gamma_n x_n+\gamma_{n-1}x_{n-1}^2+\gamma_{n-2}x_{n-2}^4+...+x_1^{2^{n-1}}\right)^{\frac{1}{2^{n-1}}}
\end{eqnarray*}
First we assume that the initial values $X_0$ of the process $X_t$ are
distributed by a smooth function $\widetilde{\rho}_0:G\longrightarrow\mathbb{R}^{\geq 0}$ which
satisfies the asymptotics 
\begin{eqnarray*}
0\leq\widetilde{\rho}_0(x_1,...,x_n)\leq C_1\psi(x_1,...,x_n)
\end{eqnarray*}
and furthermore
\begin{eqnarray*}
\int_G\frac{\widetilde{\rho}_0(x_1,...,x_n)}{\psi(x_1,...,x_n)}dx\leq C_1
\end{eqnarray*}
Let
\begin{eqnarray*}
\widetilde{\rho}(t,y_1,...,y_n)=\int_G\widetilde{\rho}_0(x_1,...,x_n)\rho(t,x_1,...,x_n,y_1,...,y_n)dx
\end{eqnarray*}
be the corresponding smooth transition density.  
\begin{theorem}
There is a constant $C_2>0$ such that for all $t\geq 0$ and all $y\in G$
\begin{eqnarray*}
0\leq\widetilde{\rho}(t,y)\leq C_2\cdot\psi(y)
\end{eqnarray*}
\end{theorem}
\textbf{Remark:} The assumption on the smooth initial distribution will be removed later on and we get an estimate
\begin{eqnarray*}
0\leq\rho(t,x_0,y)\leq C_2(x_0)\cdot\psi(y)\mbox{ for all }t\geq 1
\end{eqnarray*}
\\  
\textbf{Proof:} The Kolmogorov equation for $\rho$
\begin{eqnarray*}
\partial_t\rho=L^{*}(\rho)
\end{eqnarray*}
translates due to linearity to
\begin{eqnarray*}
\partial_t\widetilde{\rho}=L^{*}(\widetilde{\rho})
\end{eqnarray*}
with initial values $\widetilde{\rho}(0,y)=\widetilde{\rho}_0(y)$. Let
\begin{eqnarray*}
h=h(t,y)=\frac{\widetilde{\rho}(t,y_1,...,y_n)}{\psi(y_1,...,y_n)}
\end{eqnarray*}
which is a smooth function and obeys
\begin{eqnarray*}
\partial_t h=\sum_{i,j}a_{ij}\partial_i\partial_j h-
\sum_i\Big(b_i-2\sum_j a_{ij}\partial_j\log\psi\Big)\partial_i h+\frac{1}{\psi}L^{*}(\psi)\cdot h
\end{eqnarray*}
The time-uniform estimate
\begin{eqnarray*}
E\left(\frac{1}{\psi(X_t)}|X_0\right)<e^{-t}\cdot\frac{1}{\psi(X_0)}+C
\end{eqnarray*}
translates to
\begin{eqnarray*}
\int_G h(t,y)dy<C_1
\end{eqnarray*}
with a uniform constant $C_1>0$ independent of $t,T$.\\
We want to apply Moser iteration to get time-independent estimates of higher $L^m$-norms of $h$. 
For a constant $C>1$ let $h_C(t,y)=\min(C,h(t,y))$ and we first consider a
fixed time intervall $t\in[0,T]$. By choosing $C$ sufficiently large we may assume from the beginning that
$h_C(t,y)=h(t,y)$ for all $y\in K$ ($K\subset\subset G$ from condition (3)), all $t\in[0,T]$ and 
furthermore that $h_C(0,y)=h(0,y)$ for all $y\in G$ because of
the initial condition $\widetilde{\rho}(0,y)\leq C_1\psi(y)$.\\ 
For $R>2$ ($R$ is here a different
parameter than in condition (1)) let $\eta_R:\mathbb{R}^{+}\longrightarrow[0,1]$ be a smooth
function which has compact support and fulfills
\begin{eqnarray*}
\eta_R(y)=1\mbox{ for }\frac{1}{R}\leq y\leq R\mbox{ and }\eta_R(y)=0\mbox{ for }y<\frac{1}{2R}\mbox{ or }2R<y\\
\left|\frac{d\eta_R}{dy}(y)\right|\leq 4R\mbox{ for }\frac{1}{2R}<y<\frac{1}{R}\mbox{ and }
\left|\frac{d\eta_R}{dy}(y)\right|\leq\frac{4}{R}\mbox{ for }R<y<2R\\
\left|\frac{d^2\eta_R}{dy^2}(y)\right|\leq 8R^2\mbox{ for }\frac{1}{2R}<y<\frac{1}{R}\mbox{ and }
\left|\frac{d^2\eta_R}{dy^2}(y)\right|\leq\frac{8}{R}\mbox{ for }R<y<2R
\end{eqnarray*}
By abuse of notation we set $\eta_R(y_1,...,y_n)=\eta_R(y_1)\cdot...\cdot\eta_R(y_n)$. We multiply the
linear equation for $h$ above by $\eta_R(y)\cdot h_C^m(t,y)$ and integrate over $G$. Because $\eta_R(y)$ has
compact support in $G$ we may perform partial integration to obtain
\begin{eqnarray*}
\int_G a_{ij}(y)(\partial_i\partial_j h)(t,y)\cdot h_C^m(t,y)\eta_R(y)dy=\\
-m\int_G a_{ij}(y)\eta_R(y)(\partial_j h)(t,y)(\partial_i h_C)(t,y)h_C^{m-1}(t,y)dy-\\
\int_G(\partial_j h)(t,y)h_C^m(t,y)\partial_i(a_{ij}\eta_R)(y)dy=\\
-m\int_G a_{ij}(y)\eta_R(y)(\partial_j h_C)(t,y)(\partial_i h_C)(t,y)h_C^{m-1}(t,y)dy+\\
m\int_G h(t,y)(\partial_j h_C)(t,y)h_C^{m-1}(t,y)\partial_i(a_{ij}\eta_R)(y)dy+\\
\int_G h(t,y)h_C^m(t,y)\partial_i\partial_j(a_{ij}\eta_R)(y)dy=\\
=-\frac{4m}{(m+1)^2}\int_G a_{ij}(y)\eta_R(y)(\partial_j h_C^{\frac{m+1}{2}})(t,y)(\partial_i h_C^{\frac{m+1}{2}})(t,y)dy+\\
m\int_G (\partial_j h_C)(t,y)h_C^{m}(t,y)\partial_i(a_{ij}\eta_R)(y)dy+\\
\int_G h(t,y)h_C^m(t,y)\partial_i\partial_j(a_{ij}\eta_R)(y)dy=\\
-\frac{4m}{(m+1)^2}\int_G a_{ij}(y)\eta_R(y)(\partial_j h_C^{\frac{m+1}{2}})(t,y)(\partial_i h_C^{\frac{m+1}{2}})(t,y)dy+\\
\int_G h_C^m(t,y)\Big(h(t,y)-\frac{m}{m+1}h_C(t,y)\Big)\partial_i\partial_j(a_{ij}\eta_R)(y)dy
\end{eqnarray*}
and similarly
\begin{eqnarray*}
-\int_G\Big(b_i(y)-2\sum_j a_{ij}(y)\partial_j\log\psi(y)\Big)(\partial_i h)(t,y)h_C^m(t,y)\eta_R(y)dy=\\
\int_G\partial_i\Big(\eta_R\Big(b_i-2\sum_j a_{ij}\partial_j\log\psi\Big)\Big)(y)
h_C^m(t,y)\Big(h(t,y)-\frac{m}{m+1}h_C(t,y)\Big)dy
\end{eqnarray*}
So we come up with
\begin{eqnarray*}
\int_G(\partial_t h)(t,y)h_C^m(t,y)\eta_R(y)dy=\\
-\frac{4m}{(m+1)^2}\int_G\sum_{i,j}a_{ij}(y)\eta_R(y)(\partial_j h_C^{\frac{m+1}{2}})(t,y)
(\partial_i h_C^{\frac{m+1}{2}})(t,y)dy+\\
\int_G h_C^m(t,y)\Big(h(t,y)-\frac{m}{m+1}h_C(t,y)\Big)\cdot\\
\Big(\sum_{i,j}\partial_i\partial_j(a_{ij}\eta_R)+
\sum_i\partial_i(\eta_R(b_i-2\sum_j a_{ij}\partial_j\log\psi))\Big)(y)dy+\\
\int_G\Big(\frac{1}{\psi}L^{*}(\psi)\Big)(y)h(t,y)h_C^m(t,y)\eta_R(y)dy
\end{eqnarray*}
We collect all terms involving derivatives of $\eta_R$ of first or second order in a summand
$\varepsilon_{R,C,m,t}$, integrate over time and get
\begin{eqnarray*}
\int_0^T\int_G(\partial_t h)(t,y)h_C^m(t,y)\eta_R(y)dydt+
\frac{1}{m+1}\int_G h^{m+1}(0,y)\eta_R(y)dy=\\
-\frac{4m}{(m+1)^2}\int_0^T\int_G\sum_{i,j}a_{ij}(y)(\partial_j h_C^{\frac{m+1}{2}})(t,y)
(\partial_i h_C^{\frac{m+1}{2}})(t,y)\eta_R(y)dydt+\\
\int_0^T\int_G h_C^m(t,y)\Big(h(t,y)-\frac{m}{m+1}h_C(t,y)\Big)\eta_R(y)\cdot\\
\Big(\sum_{i,j}\partial_i\partial_j a_{ij}+
\sum_i\partial_i(b_i-2\sum_j a_{ij}\partial_j\log\psi)\Big)(y)dydt+\\
\int_0^T\varepsilon_{R,C,m,t}dt+
\int_0^T\int_G\Big(\frac{1}{\psi}L^{*}(\psi)\Big)(y)h(t,y)h_C^m(t,y)\eta_R(y)dydt+\\
\frac{1}{m+1}\int_G h^{m+1}(0,y)\eta_R(y)dy
\end{eqnarray*} 
We claim that for fixed $C>>1$ we have uniformly in $t\in[0,T]$
\begin{eqnarray*}
\lim_{R\rightarrow\infty}\varepsilon_{R,C,m,t}=0
\end{eqnarray*}
Consider for example ($i\neq j$)
\begin{eqnarray*}
\left|\int_G h_C^m(t,y)\Big(h(t,y)-\frac{m}{m+1}h_C(t,y)\Big)a_{ij}(y)(\partial_i\partial_j\eta_R)(y)dy\right|<\\
16C^m\int_{\mathbb{R}_{+}^{n-2}}\int_{y_i,y_j\in\mathbb{R}_{+}\backslash[\frac{1}{R},R]}
h(t,y)|a_{ij}(y)|\Big(1+\frac{1}{y_i}\Big)\Big(1+\frac{1}{y_j}\Big)dy
\end{eqnarray*}
Now according to assumption (3) we have
\begin{eqnarray*}
\int_G h(t,y)|a_{ij}(y)|\Big(1+\frac{1}{y_i}\Big)\Big(1+\frac{1}{y_j}\Big)dy<C_1<\infty
\end{eqnarray*}
and so
\begin{eqnarray*}
\lim_{R\rightarrow\infty}
16C^m\int_{\mathbb{R}_{+}^{n-2}}\int_{y_i,y_j\in\mathbb{R}_{+}\backslash[\frac{1}{R},R]}
h(t,y)|a_{ij}(y)|\Big(1+\frac{1}{y_i}\Big)\Big(1+\frac{1}{y_j}\Big)dy=0
\end{eqnarray*}
So we conclude
\begin{eqnarray*}
\lim_{R\rightarrow\infty}\varepsilon_{R,C,m,t}=0=\lim_{R\rightarrow\infty}\int_0^T\varepsilon_{R,C,m,t}dt
\end{eqnarray*}
Besides the functions
\begin{eqnarray*}
R\longmapsto-\frac{4m}{(m+1)^2}\int_G\sum_{i,j}a_{ij}(y)(\partial_j h_C^{\frac{m+1}{2}})(t,y)
(\partial_i h_C^{\frac{m+1}{2}})(t,y)\eta_R(y)dy
\end{eqnarray*}
\begin{eqnarray*}
R\longmapsto\int_G\Big(\frac{1}{\psi}L^{*}(\psi)\Big)(y)h(t,y)h_C^m(t,y)\eta_R(y)dy
\end{eqnarray*}
are decreasing for $R>R_0$ as $\frac{1}{\psi}L^{*}(\psi)<-1$ outside a fixed compact cube.
According to assumption (2) we had
\begin{eqnarray*}
\sum_{i,j}\partial_i\partial_j(a_{ij})+\sum_i\partial_i(b_i-2\sum_j a_{ij}\partial_j\log\psi)\leq 0
\end{eqnarray*}
outside of $K$ and so
\begin{eqnarray*}
R\longmapsto
\int_G h_C^m(t,y)\Big(h(t,y)-\frac{m}{m+1}h_C(t,y)\Big)\eta_R(y)\cdot\\
\Big(\sum_{i,j}\partial_i\partial_j(a_{ij})+
\sum_i\partial_i(b_i-2\sum_j a_{ij}\partial_j\log\psi)\Big)(y)dy
\end{eqnarray*}
is decreasing as well for $R>R_0$. Because of
\begin{eqnarray*}
\int_G h(t,y)dy<C_1\mbox{ uniformly in }t
\end{eqnarray*}
and $0\leq h(0,y)<C_1$ by construction we have
\begin{eqnarray*}
R\longmapsto\frac{1}{m+1}\int_G h^{m+1}(0,y)\eta_R(y)dy
\end{eqnarray*}
is increasing but bounded above by
\begin{eqnarray*}
\frac{1}{m+1}\int_G h^{m+1}(0,y)\eta_R(y)dy\leq
\frac{1}{m+1}C_1^m\int_G h(0,y)dy\leq\frac{1}{m+1}C_1^{m+1}
\end{eqnarray*}
In the next step we prove that for every $y\in G$
\begin{eqnarray*}
(R,C)\longmapsto\eta_R(y)\cdot\left(\int_0^T(\partial_t h)(t,y)h_C^m(t,y)dt+\frac{1}{m+1}h^{m+1}(0,y)\right)
\end{eqnarray*}
is an increasing function of both arguments. We know that $h(0,y)=h_C(0,y)$ for all $y\in G$.
We calculate
\begin{eqnarray*}
\int_0^T(\partial_t h)(t,y)h_C^m(t,y)\eta_R(y)dt=h(T,y)h_C^m(T,y)-h(0,y)h_C^m(0,y)-\\
m\int_0^T h(t,y)\partial_t h_C(t,y)h_C^{m-1}(t,y)dt=\\
h(T,y)h_C^m(T,y)-h^{m+1}(0,y)-m\int_0^T \partial_t h_C(t,y)h_C^{m}(t,y)dt=\\
h(T,y)h_C^m(T,y)-h^{m+1}(0,y)-\frac{m}{m+1}\left(h_C^{m+1}(T,y)-h_C^{m+1}(0,y)\right)=\\
h(T,y)h_C^m(T,y)-\frac{m}{m+1}h_C^{m+1}(T,y)-\frac{1}{m+1}h^{m+1}(0,y)
\end{eqnarray*}
and so
\begin{eqnarray*}
\int_0^T(\partial_t h)(t,y)h_C^m(t,y)dt+\frac{1}{m+1}h^{m+1}(0,y)=h(T,y)h_C^m(T,y)-\frac{m}{m+1}h_C^{m+1}(T,y)\geq 0
\end{eqnarray*}
So
\begin{eqnarray*}
R\longmapsto\eta_R(y)\cdot\left(\int_0^T(\partial_t h)(t,y)h_C^m(t,y)dt+\frac{1}{m+1}h^{m+1}(0,y)\right)
\end{eqnarray*}
is increasing. Furthermore
\begin{eqnarray*}
h(T,y)h_C^m(T,y)-\frac{m}{m+1}h_C^{m+1}(T,y)=h(T,y)C^m-\frac{m}{m+1}C^{m+1}\mbox{ if }h(T,y)>C\\
h(T,y)h_C^m(T,y)-\frac{m}{m+1}h_C^{m+1}(T,y)=\frac{1}{m+1}h^{m+1}(T,y)\mbox{ if }h(T,y)\leq C
\end{eqnarray*}
and both expressions on the right side are nondecreasing functions of $C$, so
\begin{eqnarray*}
C\longmapsto\eta_R(y)\cdot\left(\int_0^T(\partial_t h)(t,y)h_C^m(t,y)dt+\frac{1}{m+1}h^{m+1}(0,y)\right)
\end{eqnarray*}
is nondecreasing.\\
\\
Now going back to our integrated equation over time and space we see that the left side is increasing in $R$
whereas the right side contains either summands decreasing in $R$ or summands converging for $R\longrightarrow\infty$.
So we may pass to the limit $R\longrightarrow\infty$ and obtain
\begin{eqnarray*}
\int_G\Big(h(T,y)h_C^m(T,y)-\frac{m}{m+1}h_C^{m+1}(T,y)\Big)dy=\\
-\frac{4m}{(m+1)^2}\int_0^T\int_G\sum_{i,j}a_{ij}(y)(\partial_j h_C^{\frac{m+1}{2}})(t,y)
(\partial_i h_C^{\frac{m+1}{2}})(t,y)dydt+\\
\int_0^T\int_G h_C^m(t,y)\Big(h(t,y)-\frac{m}{m+1}h_C(t,y)\Big)\cdot\\
\Big(\sum_{i,j}\partial_i\partial_j a_{ij}+
\sum_i\partial_i(b_i-2\sum_j a_{ij}\partial_j\log\psi)\Big)(y)dydt+\\
\int_0^T\int_G\Big(\frac{1}{\psi}L^{*}(\psi)\Big)(y)h(t,y)h_C^m(t,y)dydt+\\
\frac{1}{m+1}\int_G h^{m+1}(0,y)dy
\end{eqnarray*} 
Now by assumption (1) there is a compact cube $K\subset\widetilde{K}\subset\subset G$ such that for
$y\in G\backslash\widetilde{K}$
\begin{eqnarray*}
\frac{1}{\psi(y)}L^{*}(\psi)(y)<-1
\end{eqnarray*}
This implies
\begin{eqnarray*}
C\longmapsto\int_G\Big(\frac{1}{\psi}L^{*}(\psi)\Big)(y)h(t,y)h_C^m(t,y)dy
\end{eqnarray*}
is decreasing in $C$ for all sufficiently large $C$.\\
So we see that the right side of the above equation is decreasing in $C$ whereas the left side is increasing in $C$.
So we may pass to the limit $C\longrightarrow\infty$ and get
\begin{eqnarray*}
\int_G h^{m+1}(T,y)dy=
-\frac{4m}{m+1}\int_0^T\int_G\sum_{i,j}a_{ij}(y)(\partial_j h^{\frac{m+1}{2}})(t,y)
(\partial_i h^{\frac{m+1}{2}})(t,y)dydt+\\
\int_0^T\int_G h^{m+1}(t,y)\cdot
\Big(\sum_{i,j}\partial_i\partial_j a_{ij}+
\sum_i\partial_i(b_i-2\sum_j a_{ij}\partial_j\log\psi)\Big)(y)dydt+\\
(m+1)\int_0^T\int_G\Big(\frac{1}{\psi}L^{*}(\psi)\Big)(y)h^{m+1}(t,y)dydt+
\int_G h^{m+1}(0,y)dy
\end{eqnarray*}
As this holds for all $T$ we infer
\begin{eqnarray*}
\frac{\partial}{\partial t}\int_G h^{m+1}(t,y)dy=
-\frac{4m}{m+1}\int_G\sum_{i,j}a_{ij}(y)(\partial_j h^{\frac{m+1}{2}})(t,y)
(\partial_i h^{\frac{m+1}{2}})(t,y)dy+\\
\int_G h^{m+1}(t,y)\cdot
\Big(\sum_{i,j}\partial_i\partial_j a_{ij}+
\sum_i\partial_i(b_i-2\sum_j a_{ij}\partial_j\log\psi)\Big)(y)dy+\\
(m+1)\int_G\Big(\frac{1}{\psi}L^{*}(\psi)\Big)(y)h^{m+1}(t,y)dy
\end{eqnarray*}
Now $(a_{ij}(y))>0$ is positive definite in every point and so
\begin{eqnarray*}
\int_G\sum_{i,j}a_{ij}(y)(\partial_j h^{\frac{m+1}{2}})(y)(\partial_i h^{\frac{m+1}{2}})(y)dy\geq
\delta_0\int_{\widetilde{K}}\left\|\nabla\left(h^{\frac{m+1}{2}}\right)\right\|^2(y)dy
\end{eqnarray*}
Furthermore
\begin{eqnarray*}
\sum_{i,j}\partial_i\partial_j a_{ij}+
\sum_i\partial_i(b_i-2\sum_j a_{ij}\partial_j\log\psi)\leq 0\mbox{ on }G\backslash\widetilde{K}
\end{eqnarray*}
and so
\begin{eqnarray*}
\frac{\partial}{\partial t}\int_G h^{m+1}(t,y)dy\leq
-\delta_0\int_{\widetilde{K}}\left\|\nabla\left(h^{\frac{m+1}{2}}\right)\right\|^2(y)dy+\\
\widetilde{C}_1(m+1)\int_{\widetilde{K}}h^{m+1}(t,y)dy-
(m+1)\int_{G}h^{m+1}(t,y)dy
\end{eqnarray*}
Here $\widetilde{C}_1>0$ is a constant independent of $T$ which can be explicitly calculated from the coefficients of the generator.
In the following we replace $\widetilde{K}$ by $K$ and $\widetilde{C}_1$ by $C_1$ for abbreviation.\\
\\
Next we use Poincare-inequality to estimate the $L^{m+1}$-norm over $K$ by the corresponding
gradient norm. Let $K_1\subset\mathbb{R}^n$ be the unit cube. Then there is a 
constant $C_n>0$ such that for any $C^1$-function
$f:K_1\longrightarrow\mathbb{R}$
\begin{eqnarray*}
\int_{K_1}\left(f-\int_{K_1}f\right)^2\leq C_n\int_{K_1}\|\nabla f\|^2
\end{eqnarray*}
Let $K_r$ be a cube with length $r$, then by rescaling we get for a $C^1$-function
$f:K_r\longrightarrow\mathbb{R}$
\begin{eqnarray*}
\int_{K_r}\left(f-\frac{1}{vol(K_r)}\int_{K_r}f\right)^2\leq C_n r^2\int_{K_r}\|\nabla f\|^2
\end{eqnarray*}
respectively
\begin{eqnarray*}
\int_{K_r}f^2\leq C_n r^2\int_{K_r}\|\nabla f\|^2+\frac{1}{r^n}\left(\int_{K_r}f\right)^2
\end{eqnarray*}
We decompose our compact cube $K\subset\subset G$ into $N=N(r)$ small cubes $K_{r,j}$ with length $r>0$ and 
find for a nonnegative $C^1$-function $f:K\longrightarrow\mathbb{R}^{\geq 0}$
\begin{eqnarray*}
\int_K f^2=\sum_{j=1}^N\int_{K_{r,j}}f^2\leq
C_n r^2\sum_{j=1}^N\int_{K_{r,j}}\|\nabla f\|^2+
\sum_{j=1}^N\frac{1}{r^n}\left(\int_{K_{r,j}}f\right)^2\leq\\
C_n r^2\int_{K}\|\nabla f\|^2+\frac{1}{r^n}\left(\int_{K}f\right)^2
\end{eqnarray*}
If we choose $r$ such that $C_n r^2=\frac{\delta_0}{C_1(m+1)}$ and $f(y)=h^{\frac{m+1}{2}}(t,y)$ we can use this to estimate
\begin{eqnarray*}
\frac{\partial}{\partial t}\int_G h^{m+1}(t,y)dy\leq-(m+1)\int_{G}h^{m+1}(t,y)dy+\\
C_1\left(\frac{C_1 C_n}{\delta_0}\right)^{\frac{n}{2}}(m+1)^{1+\frac{n}{2}}
\left(\int_K h^{\frac{m+1}{2}}(t,y)dy\right)^2\leq\\
-(m+1)\int_{G}h^{m+1}(t,y)dy+\\
C_1\left(\frac{C_1 C_n}{\delta_0}\right)^{\frac{n}{2}}(m+1)^{1+\frac{n}{2}}
\left(\int_G h^{\frac{m+1}{2}}(t,y)dy\right)^2
\end{eqnarray*}
So if we define for a moment
\begin{eqnarray*}
F(t)=\int_G h^{m+1}(t,y)dy\mbox{ and }G(t)=
C_1\left(\frac{C_1 C_n}{\delta_0}\right)^{\frac{n}{2}}(m+1)^{1+\frac{n}{2}}
\left(\int_G h^{\frac{m+1}{2}}(t,y)dy\right)^2
\end{eqnarray*}
we have
\begin{eqnarray*}
\frac{dF}{dt}\leq-(m+1)F(t)+G(t)
\end{eqnarray*}
From this we infer by Gronwall
\begin{eqnarray*}
F(t)\leq e^{-(m+1)t}\left(F(0)+\int_0^t e^{(m+1)s}G(s)ds\right)\leq\\
e^{-(m+1)t}F(0)+\Big(\sup_{t\in[0,T]}G(t)\Big)\frac{1}{m+1}(1-e^{-(m+1)t})
\end{eqnarray*}
or
\begin{eqnarray*}
\sup_{t\in[0,T]}F(t)\leq\frac{1}{m+1}\sup_{t\in[0,T]}G(t)+
\sup_{t\in[0,T]}\left(e^{-(m+1)t}\left(F(0)-\frac{1}{m+1}\sup_{t\in[0,T]}G(t)\right)\right)
\end{eqnarray*}
We may assume that the second term on the right is negative because otherwise we would have
trivial bounds of the $L^m$-norms of $h(t,y)$ only depending on the smooth bounded initial distribution.
So we may estimate
\begin{eqnarray*}
\sup_{t\in[0,T]}\int_G h^{m+1}(t,y)dy\leq
C_1^2\left(\frac{C_1^2 C_n}{2\delta_0}\right)^{\frac{n}{2}}(m+1)^{\frac{n}{2}}
\sup_{t\in[0,T]}\left(\int_G h^{\frac{m+1}{2}}(t,y)dy\right)^2
\end{eqnarray*}
Now we use that
\begin{eqnarray*}
\sup_{t\in[0,T]}\int_G h(t,y)dy\leq C_1
\end{eqnarray*}
is uniformly bounded independent of $T$.
Apply Moser iteration to end up with
\begin{eqnarray*}
\left(\sup_{t\in[0,T]}\int_G h^{2^k}(t,y)dy\right)^{2^{-k}}\leq
C_1\prod_{j=1}^k\left(C_1\left(\frac{C_1 C_n}{\delta_0}\right)^{\frac{n}{2}}2^{\frac{jn}{2}}
\right)^{2^{-j}}<C_2
\end{eqnarray*}
and so especially
\begin{eqnarray*}
\sup_{t\in[0,T]}\sup_{y\in G}h(t,y)\leq C_2
\end{eqnarray*}
is bounded independent of $T>0$. This proves the Theorem for a smooth, fast decaying initial distribution.\\
\\
In the next step we want to get rid of the assumption of a smooth initial distribution:
\begin{theorem} 
There is a constant $C_2(x_0)>0$ depending only on the initial starting point $x_0$ such that for all $t\geq 1$ 
and all $y\in G$
\begin{eqnarray*}
0\leq\rho(t,x_0,y)\leq C_2(x_0)\cdot\psi(y)
\end{eqnarray*}
\end{theorem}
\textbf{Proof:} Let $K\subset\subset G$ be a compact cube such that $\frac{1}{\psi}L^{*}(\psi)\leq -1$ on
$G\backslash K$. By enlarging $K$ if necessary we can assume that $B_2(x_0)\subset K$. Let
$\widetilde{K}=K\backslash B_1(x_0)$ and consider a family of smooth initial distribution
$\widetilde{\rho}_{\varepsilon}(t=0,y)$ with support in $B_1(x_0)$ and converging to a Dirac distribution in $x_0$
for $\varepsilon\rightarrow 0$. Especially we have
\begin{eqnarray*}
\int_G\widetilde{\rho}_{\varepsilon}(t=0,y)dy=1\mbox{ for all }\varepsilon>0
\end{eqnarray*}
Let $\widetilde{\rho}_{\varepsilon}(t,y)$ be the corresponding solution of the parabolic Kolmogorov-equation.
As before we define
\begin{eqnarray*}
0\leq h_{\varepsilon}(t,y)=\frac{\widetilde{\rho}_{\varepsilon}(t,y)}{\psi(y)}
\end{eqnarray*}
They all solve the same parabolic PDE
\begin{eqnarray*}
\partial_t h_{\varepsilon}=\widehat{L}(h_{\varepsilon})+c(y)\cdot h_{\varepsilon}
\end{eqnarray*}
where $\widehat{L}$ is a linear elliptic operator (degenerate towards the boundary $\partial G$) only involving
second and first derivatives and we have $c(y)\leq -1$ on $G\backslash K$. We know
\begin{eqnarray*}
0<\int_G h_{\varepsilon}(t,y)dy\leq C_1
\end{eqnarray*}
with a uniform constant $C_1>0$ for all $t\geq 0$ and independent of $\varepsilon>0$. This implies especially
\begin{eqnarray*}
\inf_{y\in\widetilde{K}}h_{\varepsilon}(t,y)\leq C_1
\end{eqnarray*}
uniform for all $t\geq 0$ and all $\varepsilon>0$. Because $\widehat{L}$ is uniformly elliptic on $K$ the
parabolic Harnack-inequality (see Appendix) together with the initial conditions $h_{\varepsilon}(0,y)=0$ on
$\widetilde{K}$ implies
\begin{eqnarray*}
\sup_{y\in\widetilde{K}}h_{\varepsilon}(t,y)\leq C_2\cdot
\inf_{y\in\widetilde{K}}h_{\varepsilon}(2,y)\leq C_1 C_2
\end{eqnarray*}
for all $0\leq t\leq 1$ and a uniform constant $C_2>0$ only depending on $(a_{ij}(y))$, $b_j(x)$, $c(x)$ and
not depending on $\varepsilon>0$.\\
Now we choose for comparison another initial density $\widehat{\rho}_0(y)=C_3\psi^2(y)>0$ so that
\begin{eqnarray*}
\int_G\frac{\widehat{\rho}_0(y)}{\psi(y)}dy<C_1 C_3\mbox{ and }
\sup_{y\in G}\frac{\widehat{\rho}_0(y)}{\psi(y)}dy<C_1 C_3
\end{eqnarray*}
Let $\widehat{\rho}(t,y)$ be the solution of the Kolmogorov-equation
with $\widehat{\rho}(0,y)=\widehat{\rho}_0(y)$ and let $\widehat{h}(t,y)=\widehat{\rho}(t,y)/\psi(y)$.
Then by Harnack-inequality we have
\begin{eqnarray*}
\inf_{y\in K,0\leq t\leq 1}\widehat{h}(t,y)>\delta_3 C_3>0
\end{eqnarray*}
Choose $C_3>0$ so large that $\delta_3 C_3>C_1 C_2$. We note that $C_3$ is still independent of $\varepsilon>0$. 
Then we claim that
\begin{eqnarray*}
0\leq h_{\varepsilon}(t,y)\leq\widehat{h}(t,y)
\end{eqnarray*}
for all $0\leq t\leq 1$, for all $y\in G\backslash K$ and all $\varepsilon>0$.\\
By construction we have $(h_{\varepsilon}-\widehat{h})(0,y)<0$ for all $y\in G\backslash K$ and
$(h_{\varepsilon}-\widehat{h})(t,y)<0$ for all $(t,y)\in[0,1]\times\partial K$. 
Furthermore by Theorem 1 we have
\begin{eqnarray*}
\lim_{y\rightarrow\partial G}(h_{\varepsilon}-\widehat{h})(t,y)=0
\end{eqnarray*}
for all $t\geq 0$ and
\begin{eqnarray*}
\partial_t(h_{\varepsilon}-\widehat{h})=\widehat{L}(h_{\varepsilon}-\widehat{h})+
c(y)\cdot(h_{\varepsilon}-\widehat{h})
\end{eqnarray*} 
Assume by contradiction that for some $0<t\leq 1$ there are $y\in G\backslash K$ with
$(h_{\varepsilon}-\widehat{h})(t,y)>0$. The set of those $t$ is open in $[0,1]$ and for those $t$ the map
$y\in G\backslash K\longmapsto(h_{\varepsilon}-\widehat{h})(t,y)$ has a strictly positive maximum in a
point $y_t\in G\backslash K$. As $(h_{\varepsilon}-\widehat{h})(0,y)<0$ and $(h_{\varepsilon}-\widehat{h})(t,y_t)>0$
there must be points $(t,y_t)$ with $\frac{d}{dt}(h_{\varepsilon}-\widehat{h})(t,y_t))>0$. But the parabolic
Kolmogorov-equation implies
\begin{eqnarray*}
\frac{d}{dt}(h_{\varepsilon}-\widehat{h})(t,y_t)\leq c(y_t)\cdot(h_{\varepsilon}-\widehat{h})(t,y_t)<0
\end{eqnarray*}
a contradiction and so we have
\begin{eqnarray*}
0\leq h_{\varepsilon}(t,y)\leq\widehat{h}(t,y)<C_5\psi(y)
\end{eqnarray*}
for all $0\leq t\leq 1$, all $y\in G\backslash K$ and all $\varepsilon>0$.\\
For estimation inside of $K$ and especially on $B_1(x_0)$ one uses a smooth bumping function 
$\eta:B_2(x_0)\longrightarrow[0,1]$ and gets from the parabolic PDE
\begin{eqnarray*}
\partial_t\int_K\eta^{2m}h_{\varepsilon}^{2m}\leq-\delta_1\int_K\left|\nabla\left(
\eta^m h_{\varepsilon}^m\right)\right|^2+2mC_1\int_K h_{\varepsilon}^{2m}
\end{eqnarray*}
The parabolic Harnack-inequality on $\widetilde{K}=K\backslash B_1(x_0)$ together with the uniform $L^1$-bounds
allows us to estimate $h_{\varepsilon}(t,y)<C_0$ for all $y\in\widetilde{K}$, all $t\in[0,1]$ and all
$\varepsilon>0$. We may assume that
\begin{eqnarray*}
\int_K\eta^m h_{\varepsilon}^m\geq C_0^m
\end{eqnarray*}
because otherwise we would have trivial bounds of the $L^m$-norm on $B_1(x_0)$. This implies
\begin{eqnarray*}
\int_K h_{\varepsilon}^m\leq C_2\int_K\eta^m h_{\varepsilon}^m
\end{eqnarray*} 
Furthermore by Nash-inequality with $\beta=\frac{2}{n}>0$
\begin{eqnarray*}
\left(\int_K\eta^{2m}h_{\varepsilon}^{2m}\right)^{1+\beta}\leq C_3\left(\int_K\left|\nabla\left(
\eta^m h_{\varepsilon}^m\right)\right|^2\right)\cdot\left(
\int_K\eta^m h_{\varepsilon}^m\right)^{2\beta}
\end{eqnarray*}
So we get
\begin{eqnarray*}
\partial_t\int_K\eta^{2m}h_{\varepsilon}^{2m}\leq-\frac{\delta_1}{C_3}
\left(\int_K\eta^{2m}h_{\varepsilon}^{2m}\right)^{1+\beta}\cdot\left(
\int_K\eta^m h_{\varepsilon}^m\right)^{-2\beta}
+2mC_1 C_2\int_K \eta^{2m}h_{\varepsilon}^{2m}
\end{eqnarray*}
So if we define
\begin{eqnarray*}
g_{m,\varepsilon}(t)=e^{-mC_1 C_2 t}\cdot\int_K\eta^m(y) h_{\varepsilon}^m(t,y)dy>0
\end{eqnarray*}
we have a differential inequality
\begin{eqnarray*}
\frac{dg_{2m,\varepsilon}}{dt}\leq-\frac{\delta_1}{C_3}g_{m,\varepsilon}^{-2\beta}\cdot
g_{2m,\varepsilon}^{1+\beta}
\end{eqnarray*}
from which we conclude
\begin{eqnarray*}
(g_{2m,\varepsilon}(t))^{-\beta}\geq(g_{2m,\varepsilon}(0))^{-\beta}+\frac{\delta_1\beta}{C_3}
\int_0^t(g_{m,\varepsilon}(s))^{-2\beta}ds
\end{eqnarray*}
Setting $m=2^k$, $k\in\mathbb{N}_0$ and using $g_{1,\varepsilon}(t)\leq C_1$ uniformly bounded one can prove by
induction that
\begin{eqnarray*}
(g_{m,\varepsilon}(t))^{-\beta}\geq(g_{m,\varepsilon}(0))^{-\beta}+
2^{-N_k}\Big(\frac{\delta_1\beta}{C_3}\Big)^{m-1} t^{m-1}
\end{eqnarray*}
with
\begin{eqnarray*}
N_k=2^k\cdot\left(\sum_{j=2}^k j\cdot 2^{-j}\right)
\end{eqnarray*}
Now letting $k\longrightarrow\infty$ yields
\begin{eqnarray*}
\sup_{y\in B_1(x_0)}h_{\varepsilon}(t,y)<2^{\frac{1}{\beta}\sum_2^{\infty}j\cdot 2^{-j}}\cdot
\Big(\frac{\delta_1\beta}{C_3}\Big)^{-\frac{1}{\beta}}\cdot t^{-\frac{1}{\beta}}\cdot
e^{C_1 C_2 t}
\end{eqnarray*}
This proves Theorem 2 as in $t=1$ we have a uniformly bounded transition density decaying fast towards $\partial G$
for which we can apply Theorem 1 for all $t\geq 1$.\\
\\
\small
Authors address:\\
Bert Koehler, Debeka Hauptverwaltung, Ferdinand-Sauerbruch-Str. 18, 56058 Koblenz, Germany,
Email: Bert.Koehler@debeka.de\\
\\
Volker Krafft, Debeka Hauptverwaltung, Ferdinand-Sauerbruch-Str. 18, 56058 Koblenz, Germany,
Email: Volker.Krafft@debeka.de\\
\\
\textbf{References}\\
(1) L.H\"ormander, Hypoelliptic second order differential equations, Acta Math. 119, 147-171, 1967\\
(2) L.C.Evans, Partial Differential Equations, AMS 1998\\
(3) E.B.Fabes and D.W.Stroock, A new Proof of Moser Parabolic Harnack Inequality via the old ideas of Nash,
Arch. Ratl. Mech. 1986\\
(4) D.Gilbarg and N.Trudinger, Elliptic Partial Differential Equations of Second Order, Springer

\end{document}